\documentclass[11pt]{article}
\usepackage{amsfonts}
\usepackage{amssymb}
\usepackage{epsfig}
\usepackage{amsmath}

\newcommand{\ignore}[1]{}
\ignore{The text in between these parentheses is ignored by LaTeX
ignore}

\def\@begintheorem#1#2{\par\bgroup{\sc #1\ #2. }\it\ignorespaces}
\def\@opargbegintheorem#1#2#3{\par\bgroup{\sc #1\ #2\ (#3). } \it\ignorespaces}
\def\@endtheorem{\egroup}
\newtheorem{theorem}{Theorem}[section]
\newtheorem{corollary}[theorem]{Corollary}
\newtheorem{lemma}[theorem]{Lemma}
\newtheorem{proposition}[theorem]{Proposition}
\newtheorem{problem}[theorem]{Problem}
\newtheorem{example}[theorem]{Example}
\newtheorem{algorithm}[theorem]{Algorithm}
\newtheorem{definition}[theorem]{Definition}
\newcommand{\bt}[1]{\begin{theorem}\label{#1}}
\newcommand{\bc}[1]{\begin{corollary}\label{#1}}
\newcommand{\bl}[1]{\begin{lemma}\label{#1}}
\newcommand{\bp}[1]{\begin{proposition}\label{#1}}
\newcommand{\bpro}[1]{\begin{problem}\label{#1}}
\newcommand{\be}[1]{\begin{example}\rm\label{#1}}
\newcommand{\ba}[1]{\begin{algorithm}\rm\label{#1}}
\newcommand{\bd}[1]{\begin{definition}\rm\label{#1}}

\newcommand{\et}{\end{theorem}}
\newcommand{\ec}{\end{corollary}}
\newcommand{\el}{\end{lemma}}
\newcommand{\ep}{\end{proposition}}
\newcommand{\epro}{\end{problem}}
\newcommand{\ee}{\end{example}}
\newcommand{\ea}{\end{algorithm}}
\newcommand{\ed}{\end{definition}}

\def\N{\mathbb{N}}
\def\R{\mathbb{R}}
\def\Z{\mathbb{Z}}

\def \G {{{\cal G}}}
\def \F {{{\cal F}}}
\def \L {{{\cal L}}}

\def \conv {{\rm conv}}
\def \zone {{\rm zone}}
\def \supp {{\rm supp}}

\newcommand{\boproof}{\noindent {\em Proof. }}
\newcommand{\eoproof}{\hspace*{\fill} $\square$ \vspace{5pt}}
\newcommand{\red}{\sqsubseteq}

\begin{document}

\title{\bf Convex Integer Maximization via Graver Bases}
\author{J. A. De Loera
\thanks{Supported in part by NSF grant DMS-0608785.}
\and R. Hemmecke
\thanks{Supported in part by the European TMR network ADONET 504438.}
\and S. Onn \thanks{Supported in part by a grant from
ISF - the Israel Science Foundation.}
\and U.G. Rothblum \thanks{Supported in part by a grant from
ISF - the Israel Science Foundation.}
\and R. Weismantel
\thanks{Supported in part by the European TMR network ADONET 504438.}}

\date{}
\maketitle

\begin{abstract}
We present a new algebraic algorithmic scheme to solve {\em convex integer maximization}
problems of the following form, where $c$ is a convex function
on $\R^d$ and $w_1 x,\dots,w_dx$ are linear forms on $\R^n$,
$$\max\, \{c(w_1 x,\dots,w_d x):\ Ax=b,\ x\in\N^n\}\ .$$
This method works for arbitrary input data $A,b,d,w_1,\dots,w_d,c$.
Moreover, for fixed $d$ and several important classes of programs in
{\em variable dimension}, we prove that our algorithm runs in {\em polynomial time}.
As a consequence, we obtain polynomial time algorithms for various types of multi-way
transportation problems, packing problems, and partitioning problems in variable dimension.
\vskip.2cm
\noindent
{\em keywords:} Graver basis, Gr\"obner basis, Graver complexity,
contingency table, transportation polytope, transportation problem,
integer programming, discrete optimization, packing, cutting stock,
partitioning, clustering, polyhedral combinatorics,
convex optimization, computational complexity.

\end{abstract}

\section{Introduction}
\label{Introduction}

In the past fifteen years algebraic geometry and commutative algebra
tools have shown their exciting potential to study problems in integer
optimization (see \cite{bertweis,rekhasurvey} and references therein).
But, so far, algebraic methods have always been considered ``guilty''
of bad computational complexity, namely, the notorious bad complexity
for computing general Gr\"obner bases when the number of variables
grow (see \cite{mayrsurvey} and references therein).  This paper
demonstrates that, by carefully analyzing the structure of toric ideals
in particular problems, algebraic tools can compete (and win!) against
more mainstream tools in optimization.

The main algebraic ingredient we will need is the notion of
\emph{Graver bases}, a special kind of universal Gr\"obner bases for
the toric ideals associated with integer matrices. We recommend the
introduction presented in Chapter 4 of \cite{berndbook} for a basic
introduction to Gr\"obner and Graver bases of toric ideals. We
consider a new algorithmic scheme for solving the following
far-reaching generalization of standard linear integer programming:

\vskip.2cm\noindent
{\bf Convex Integer Maximization}. Given positive integers
$d,m,n$, integer vectors $w_1,\dots,w_d\in\Z^n$ and $b\in\Z^m$, integer
matrix $A\in \Z^{m\times n}$, and convex function $c:\R^d\longrightarrow\R$,
find a nonnegative integer vector $x\in\N^n$ maximizing the objective function
$c(w_1 x,\dots,w_d x)$ subject to the equation system $Ax=b$,
$$\max\, \{c(w_1 x,\dots,w_d x):\ Ax=b,\ x\in\N^n\}\ .$$
This problem can be interpreted as multi-objective integer
programming: given $d$ different linear objective functions
$w_1,\dots,w_d$, the goal is to maximize their ``convex balancing"
given by $c(w_1 x,\dots,w_d x)$.  The convex integer maximization
problem is very expressive and in fact, contains a whole hierarchy of
problems of increasing generality and complexity, parameterized by the
number $d$ of linear objectives used: at the bottom lies the standard
linear integer programming problem, recovered as the special case of
$d=1$ and $c$ the identity on $\R$; and at the top lies the problem of
maximizing an arbitrary convex functional over the set of integer
points in a rational polyhedron in $\R^n$, arising with $d=n$ and
$w_i={\bf 1}_i$ the $i$-th standard unit vector in $\R^n$ for all $i$.

In general, convex integer maximization is intractable even for small
fixed $d$, since already for $d=1$ it includes linear integer
programming which is NP-hard.  For variable $d$, even very simple
special cases are NP-hard, such as the following instance ({\em
positive semi-definite quadratic binary programming}),
$$\max\, \{(w_1 x)^2+\cdots+ (w_n x)^2:
\ x_i+y_i=1\ \ (i=1,\dots,n),\ \ x,y\in\N^n\}\ .$$
Clearly, the complexity of the problem depends also on the presentation
of the convex function: we will assume that $c$ is presented
by a {\em comparison oracle} that, queried on $x,y\in\R^d$,
asserts whether or not $c(x)\leq c(y)$.
This is a very broad presentation that reveals little information
on the function, making the problem harder to solve. In particular,
if the polyhedron $\{x\in\R^n_+\,:\, Ax=b\}$ is unbounded, then the problem
is inaccessible even in one variable with no equation constraints:
consider the following family of univariate convex integer programs
with convex functions parameterized by $-\infty< u\leq\infty$,
$$
\max\, \{c_u(x)\ :\ x\in\N\}\ ,\quad
c_u(x):=\left\{
\begin{array}{ll}
    -x, & \hbox{if   $x<u$;} \\
    x-2u, & \hbox{if $x\geq u$.} \\
\end{array}
\right.\ ;
$$
now consider any algorithm attempting to solve the problem and let $u$
be the maximum value of $x$ in all queries to the oracle of $c$; then
the algorithm can not distinguish between the problem with $c_u$,
whose objective function is unbounded, and the problem with
$c_\infty$, whose optimal objective value is $0$.
(We remark that, for explicitly given (rather than oracle presented)
simple convex functions, it might be possible to handle unbounded
feasible regions as well; this should be the subject of future study.)

In spite of these difficulties, we show in this article that the
algebraic techniques of Graver bases allow us to solve the convex integer
maximization problem in polynomial time for a large and useful class
of integer programs in variable dimension. Moreover, this class is
{\em universal} for integer programming in a well defined sense,
enabling to extend this to an algorithmic scheme for solving convex
integer maximization over arbitrary integer programs.

Our first key lemma, extending results of \cite{OR} for combinatorial
optimization, shows that when a suitable geometric condition holds, it
is possible to efficiently reduce the convex integer maximization
problem to the solution of polynomially many linear integer
programming counterparts.  As we will see, this condition holds
naturally for a broad class of problems in variable dimension. To
state this result, we need the following terminology. A {\em direction}
of an edge ($1$-face) $e$ of a polyhedron $P$ is an
nonzero scalar multiple of $u-v$ with $u,v$ any two distinct points in
$e$.  A set of vectors {\em covers all edge-directions of $P$} if it
contains a direction of each edge of $P$.  A {\em linear integer
  programming oracle} for matrix $A\in\Z^{m\times n}$ and vector
$b\in\Z^m$ is one that, queried on $w\in\Z^n$, solves the linear
integer program $\max\{wx: Ax=b,\ x\in\N^n\}$, that is, either returns
an optimal solution $x\in\N^n$, or asserts that the program is
infeasible, or asserts that the objective function $w$ is unbounded.

\bl{EdgeDirections}
For any fixed $d$ there is a strongly polynomial oracle-time algorithm that,
given any vectors $w_1,\dots,w_d\in\Z^n$, matrix $A\in \Z^{m\times n}$ and
vector $b\in\Z^m$ endowed with a linear integer programming oracle,
finite set $E\subset\Z^n$ covering all edge-directions of the polyhedron
$\conv\{x\in\N^n\,:\, Ax=b\}$, and convex functional $c:\R^d\longrightarrow\R$
presented by a comparison oracle, solves the convex integer program
$$\max\, \{c(w_1 x,\dots,w_d x):\ Ax=b,\ x\in\N^n\}\ .$$
\el

Here, {\em solving} the program means that the algorithm
either returns an optimal solution $x\in\N^n$, or asserts
the problem is infeasible, or asserts the polyhedron
$\{x\in\R^n_+\,:\, Ax=b\}$ is unbounded in which case
the problem is hopeless (see discussion above); and {\em strongly polynomial
oracle-time} means that the number of arithmetic operations and calls to
the oracles are polynomially bounded in $m$ and $n$, and the size of the
numbers occurring throughout the algorithm is polynomially bounded in
the size of the input (which is the number of bits in the binary
representation of the entries of $w_1,\dots,w_d,A,b,E$).

Our main theorem, building on Lemma \ref{EdgeDirections}, shows
that a broad (in fact, {\em universal}) class of convex integer
maximization problems can be solved in polynomial time.
Given an $(r+s)\times t$ matrix $A$,
let $A_1$ be its $r\times t$ sub-matrix consisting of the first $r$
rows and let $A_2$ be its $s\times t$ sub-matrix consisting of the
last $s$ rows. Define the {\em n-fold matrix} of $A$ to be the
following $(r+ns)\times nt$ matrix,
$$A^{(n)}\quad:=\quad ({\bf 1}_n\otimes A_1)\oplus(I_n \otimes A_2)\quad=\quad
\left(
\begin{array}{ccccc}
  A_1    & A_1    & A_1    & \cdots & A_1    \\
  A_2  & 0      & 0      & \cdots & 0      \\
  0  & A_2      & 0      & \cdots & 0      \\
  \vdots & \vdots & \ddots & \vdots & \vdots \\
  0  & 0      & 0      & \cdots & A_2      \\
\end{array}
\right)\quad .
$$
Note that $A^{(n)}$ depends on $r$ and $s$: these will be indicated by
referring to $A$ as an ``$(r+s)\times t$ matrix".

We establish the following theorem, which asserts that
convex integer maximization over $n$-fold systems of a fixed matrix $A$,
in variable dimension $nt$, are solvable in polynomial time.
This extends results for {\em linear} integer programming from \cite{DHOW}.

\bt{Main}
For any fixed positive integer $d$ and fixed $(r+s)\times t$ integer
matrix $A$ there is a polynomial oracle-time algorithm that,
given $n$, vectors $w_1,\dots,w_d\in\Z^{nt}$ and $b\in\Z^{r+ns}$,
and convex function $c:\R^d\longrightarrow\R$ presented by a
comparison oracle, solves the convex n-fold integer maximization problem
$$\max\, \{c(w_1 x,\dots,w_d x):\ A^{(n)}x=b,\ x\in\N^{nt}\}\ .$$
\et

The equations defined by an n-fold matrix have the following,
perhaps more illuminating, interpretation: splitting the variable vector
and the right-hand side vector into components of suitable sizes,
$x=(x^1,\dots,x^n)$ and $b=(b^0,b^1,\dots,b^n)$, where $b^0\in\Z^r$
and $x^k\in\N^t$ and $b^k\in \Z^s$ for $k=1,\dots,n$, the equations
become $A_1(\sum_{k=1}^n x^k)=b^0$ and $A_2x^k=b^k$ for $k=1,\dots,n$.
Thus, each component $x^k$ satisfies a system of constraints defined
by $A_2$ with its own right-hand side $b^k$, and the sum
$\sum_{k=1}^n x^k$ obeys constraints determined by $A_1$ and $b^0$
restricting the ``common resources shared by all components".

Theorem \ref{Main} has various applications, including to multiway
transportation problems, packing problems, vector partitioning and
clustering, which will be discussed in Section \ref{Applications}.
For example, we have the following corollary providing the
first polynomial time solution of convex $3$-way transportation.

\bc{Threeway}{\bf (convex 3-way transportation)}
For any fixed $d,p,q$ there is a polynomial oracle-time algorithm that,
given $n$, arrays $w_1,\dots,w_d\in\Z^{p\times q\times n}$,
$u\in\Z^{p\times q}$, $v\in\Z^{p\times n}$, $z\in\Z^{q\times n}$,
and convex $c:\R^d\longrightarrow\R$ presented by comparison oracle,
solves the convex integer 3-way transportation problem
$$\max\{\,c(w_1 x, \dots, w_d x)
 \ :\ x\in\N^{p\times q\times n}\,,\ \sum_i x_{i,j,k}=z_{j,k}
\,,\ \sum_j x_{i,j,k}=v_{i,k}\,,\ \sum_k x_{i,j,k}=u_{i,j}\,\}\ .$$
\ec
Note that in contrast, if the dimensions of two sides of the tables are variable,
say, $q$ and $n$, then even the standard {\em linear} integer $3$-way
transportation problem over such tables is NP-hard, see \cite{DO1,DO2,DO3}.

\vskip.5cm

We proceed to discuss the universality of $n$-fold integer programming and
describe our algorithmic scheme for solving convex integer maximization
over an arbitrary system. Define a variant of the $n$-fold operator as follows:
for an $s\times t$ matrix $A$, define its {\em $n$-product} $A^{[n]}$ to be
the $n$-fold product of the $(t+s)\times t$ matrix obtained by
appending $A$ to the $t\times t$ identity matrix $I_t$, that is:
$$A^{[n]}\ :=\
\left(
\begin{array}{c}
  I_t  \\
  A    \\
\end{array}
\right)^{(n)}\ .
$$
Consider $m$-products $(1,1,1)^{[m]}$ of the $1\times 3$ matrix $(1,1,1)$.
Note that $(1,1,1)^{[m]}$ is precisely the $(3+m)\times 3m$ vertex-edge
incidence matrix of the complete bipartite graph $K_{3,m}$. For instance,
$$(1,1,1)^{[3]}\ =\
\left(
\begin{array}{ccccccccc}
  1 & 0 & 0 & 1 & 0 & 0 & 1 & 0 & 0 \\
  0 & 1 & 0 & 0 & 1 & 0 & 0 & 1 & 0 \\
  0 & 0 & 1 & 0 & 0 & 1 & 0 & 0 & 1 \\
  1 & 1 & 1 & 0 & 0 & 0 & 0 & 0 & 0 \\
  0 & 0 & 0 & 1 & 1 & 1 & 0 & 0 & 0 \\
  0 & 0 & 0 & 0 & 0 & 0 & 1 & 1 & 1 \\
\end{array}
\right)\ .
$$
The following result which incorporates the recent universality theory
of \cite{DO1,DO2,DO3} asserts that {\em every} convex integer maximization
problem can be lifted in polynomial time to some convex integer maximization
problem defined by some $n$-product of some $m$-product of $(1,1,1)$.
Theorem \ref{Main} can then be harnessed to solve the lifted program,
providing a general solution scheme for convex maximization.

\bc{Lifting}{\bf (scheme for arbitrary convex integer maximization)}
There is a polynomial time algorithm that, given integer $p\times q$
matrix $B$ and $b\in\Z^p$ with $\{x\in\Z^q\,:\,Bx=b\,,\ x\geq 0\}$ bounded,
computes $m$, $n$, and integer $(3m+n(3+m))$ vector $a$
such that, for any given vectors $w_1,\dots w_d\in\Z^q$
and any convex function $c$ on $\R^d$, the corresponding
convex integer maximization problem lifts to a convex integer maximization
problem defined by the $n$-product of the $m$-product of $(1,1,1)$, that is,
\begin{eqnarray*}
& & \max\{c(w_1x,\dots,w_dx)\,:\,x\in\Z^q\,,\ Bx=b\,,\ x\geq 0\} \\
& = & \max\left\{c({\hat w}_1 {\hat x},\dots,{\hat w}_d{\hat x})
\ :\ {\hat x}\in\Z^{3mn}\,,\ \ \left((1,1,1)^{[m]}\right)^{[n]}{\hat x}
\, =\, a\,,\ \ {\hat x}\geq 0\right\}\ .
\end{eqnarray*}
\ec
The algorithm also computes an embedding of $\Z^q$ into $\Z^{3mn}$ so that the vectors
${\hat w}_1,\dots,{\hat w}_d\in\Z^{3mn}$ are obtained from the corresponding
vectors $w_1,\dots w_d\in\Z^q$ by simply adding sufficiently many $0$ entries.

\boproof
Reformulating the universality theorem for multiway tables from \cite{DO2} in
terms of products, it asserts that the set of integer points
$\{x\in\Z^q\,:\, Bx=b\,,\ x\geq 0\}$ in any rational polytope
stands in polynomial-time computable coordinate-embedding linear bijection
with the set of integer points in the polytope
$\left\{{\hat x}\in\Z^{3mn}\,:\, \left((1,1,1)^{[m]}\right)^{[n]}{\hat x}
\, =\, a\,,\ {\hat x}\geq 0\right\}$ for some $m$, $n$ and $a$.
Lifting each $w_i\in\Z^q$ to ${\hat w}_i\in\Z^{3mn}$ by adding suitable
$0$ entries, implies that for every integer point $x$ in the original program
and its corresponding integer point ${\hat x}$ in the lifted program,
we have the same objective function value
$c(w_1x,\dots,w_dx)=c({\hat w}_1 {\hat x},\dots,{\hat w}_d{\hat x})$.
Thus, the optimal objective function values in the original and lifted
programs are the same, and, moreover, an optimal solution to the original
program can be read off as any point $x$ corresponding to any
optimal solution ${\hat x}$ to the lifted program.
\eoproof

Note that, if P$\neq$NP, there can be no polynomial time algorithm
for general linear integer programming, let alone convex integer maximization.
So how does this reconcile with the scheme suggested by Corollary \ref{Lifting}
above ? The point is that, for every fixed $m$, Theorem \ref{Main} provides
a polynomial time algorithm for convex maximization over all integer programs
that lift to programs with defining matrix that is the $n$-product
$\left((1,1,1)^{[m]}\right)^{[n]}$ of $(1,1,1)^{[m]}$.
But for arbitrary integer programs, $m$ is variable as well and so
the whole procedure is not polynomial. But in practice, this might
be efficient or enable a quick approximation, and should be
the subject of future study. We also note that, for fixed $m$, the computational
complexity of solving convex maximization over programs defined by
$\left((1,1,1)^{[m]}\right)^{[n]}$ is dominated by $n^{dg(m)}$,
where $g(m)$ is the so-called {\em Graver complexity} of the complete
bipartite graph $K_{3,m}$ and of its incidence matrix
$(1,1,1)^{[m]}$. The precise rate of growth of $g(m)$
as a function of $m$ is unknown and intriguing; see \cite{BO} for
the best bounds and for more details and precise definitions.

\vskip.3cm

The rest of the article proceeds as follows. In Section \ref{Proofs}
we give the proofs of all statements. We begin by discussing
edge-directions of polyhedra and provide the algorithm establishing
Lemma \ref{EdgeDirections}. We proceed to discuss Graver bases and,
incorporating Lemma \ref{EdgeDirections} and recent results from
\cite{DHOW}, which are based on results of Ho\c{s}ten and Sullivant
\cite{HS} and Santos and Sturmfels \cite{SS} on the asymptotic
stabilization of Graver bases, we are able to establish
Theorem \ref{Main}. In Section \ref{Applications} we discuss
applications to multiway transportation, packing, vector
partitioning and clustering, as follows: in
\ref{Applications1} we obtain Corollary \ref{Threeway} and an
extension to $k$-way transportation problems of any dimension $k$
(Corollary \ref{Multiway}); in \ref{Applications2} we describe
applications to bin packing problems (Corollary \ref{Packing});
finally, in \ref{Applications3} we apply our Theorem \ref{Main} to vector
partitioning in general and clustering in particular (Corollary \ref{Partition}).

\section{Proofs}
\label{Proofs}

In this section we prove Lemma \ref{EdgeDirections}, which is of interest in its own
right, and combine it with several other results to establish our main Theorem \ref{Main}.
Before proceeding with the details, we provide the main outline and point out the
difficulties that we have to overcome. Given data for a convex integer maximization
problem $\max\{c(w_1x,\dots,w_dx):\ Ax=b,\ x\in\N^n\}$, consider the polyhedron
$P:=\conv\{x\in\N^n:Ax=b\}\subseteq\R^n$ and its projection
$Q:=\{(w_1x,\dots,w_dx):\ x\in P\}\subseteq\R^d$. Note that $P$ is
the so-called {\em integer hull} of $\{x\in\R^n:Ax=b\,,\ x\geq 0\}$
and has typically exponentially many vertices and is not accessible
computationally. Note also that, since $c$ is convex, there is an optimal
solution $x$ whose projection $(w_1x,\dots,w_dx)$ is a vertex of $Q$.
So an important ingredient in the solution is to construct the vertices of $Q$.
Unfortunately, $Q$ may also have exponentially many vertices even though
it lives in a space $\R^d$ of fixed dimension. However, we will be able
to show that, when the number of {\em edge-directions}
of $P$ is polynomial, the number of vertices of $Q$ is polynomial.
Nonetheless, even in this case, it is not possible to construct these vertices
directly, since the number of vertices of $P$ may still be exponential.
To overcome this difficulty, we need to make use of a suitable {\em zonotope}.
This is the key idea underlying the algorithm of Lemma \ref{EdgeDirections}.
Next, we restrict attention to $n$-fold systems. For such systems,
using recent results of \cite{HS,SS} on the stabilization of their {\em Graver bases},
we are able to show that the set of edge-directions of the integer hull $P$
can be computed in polynomial time. Combining this with Lemma \ref{EdgeDirections}
and several other results from \cite{DHOW} we obtain Theorem \ref{Main}.

We now proceed with the precise details.
As defined earlier, a {\em direction} of an edge ($1$-face) $e$ of a
polyhedron $P$ is any nonzero scalar multiple of $u-v$ where $u,v$
are any two distinct points in $e$. We say that a set of vectors $E$
{\em covers all edge-directions of $P$} if it contains a direction of
each edge of $P$. A polyhedron $Z$ is a {\em refinement} of a polyhedron
$P$ if the closure of each normal cone of $P$ is the union of closures of
normal cones of $Z$. The {\em zonotope} generated by a finite
set $E\subset\R^n$ is the polytope
$Z:=\zone(E):=\conv\{\sum_{e\in E} \lambda_e e:\lambda_e=\pm 1\}$.
More details and proofs of the next two propositions can
be found in \cite{GS,OR,OS} and the references therein.

\bp{Refinement}
Let $E\subset\R^n$ be a finite set covering all
edge-directions of a polyhedron $P\subseteq\R^n$.
Then the {\em zonotope}
$Z:=\zone(E)=\conv\{\sum_{e\in E} \lambda_e e:\lambda_e=\pm 1\}$
generated by $E$ is a refinement of $P$.
\ep

\bp{Zonotopes}
For any fixed $d$, there is a polynomial time algorithm that, given
any $E\subset\Z^d$, outputs every vertex $v$ of $Z:=\zone(E)$ along
with $g_v\in\Z^d$ with $g_v x$ uniquely maximized over $Z$ at $v$.
\ep
We can now prove Lemma \ref{EdgeDirections},
showing that a set of edge-directions of the polyhedron underlying a
convex integer program allows to solve it by solving
polynomially many linear integer counterparts.

\vskip.2cm\noindent{\bf Lemma \ref{EdgeDirections}}
{\em For any fixed $d$ there is a strongly polynomial oracle-time algorithm that,
given any \break vectors $w_1,\dots,w_d\in\Z^n$, matrix $A\in \Z^{m\times n}$
and vector $b\in\Z^m$ endowed with a linear integer programming oracle,
finite set $E\subset\Z^n$ covering all edge-directions of the polyhedron
$\conv\{x\in\N^n\,:\, Ax=b\}$, and convex functional $c:\R^d\longrightarrow\R$
presented by a comparison oracle, solves the convex integer program
$$\max\, \{c(w_1 x,\dots,w_d x):\ Ax=b,\ x\in\N^n\}\ .$$}

\boproof
We provide the algorithm claimed by the theorem.
First, query the linear integer programming oracle of $A,b$ on the trivial
linear function $w=0$; if the oracle asserts that the linear problem is
infeasible, then terminate the algorithm asserting that the convex problem
is infeasible. So assume the
problem is feasible. Let $P:=\conv\{x\in\N^n:Ax=b\}\subseteq\R^n$
and $Q:=\{(w_1x,\dots,w_dx):\ x\in P\}\subseteq\R^d$.
Then $Q$ is a projection of $P$, and the corresponding projection
$D:=\{(w_1e,\dots,w_de):\ e\in E\}$ of the set $E$ is a set covering all
edge-directions of $Q$. Let $Z:=\zone(D)\subset\R^d$ be the zonotope
generated by $D$. Since $d$ is fixed, by Proposition \ref{Zonotopes}
we can produce in polynomial time all vertices of $Z$, every vertex
$v$ along with $g_v\in\Z^d$ such that the linear function defined by
$g_v$ is uniquely maximized over $Z$ at $v$. For each of the polynomially
many $g_v$, repeat the following procedure. Define a vector $h_v\in\Z^n$
by $h_{v,j}:=\sum_{i=1}^d w_{i,j}g_{v,i}$ for $j=1,\dots,n$.
Now query the linear integer programming oracle of $A,b$ on the linear
function $w:=h_v\in\Z^n$. If the oracle replies that the objective is
unbounded, then terminate the algorithm asserting that $P$ is an unbounded
polyhedron. Otherwise, let $x_v\in P\cap\N^n$ be the optimal solution obtained
from the oracle, and let $z_v:=(w_1x_v,\dots,w_dx_v)\in Q$ be its projection.
Since for every $x\in P$ and its projection $z:=(w_1x,\dots,w_dx)\in Q$
we have $g_vz=h_vx$, we conclude that $z_v$ is a maximizer of $g_v$ over $Q$.
Now we claim that each vertex $u$ of $Q$ equals some $z_v$. Indeed,
since $Z$ is a refinement of $Q$ by Proposition \ref{Refinement},
it follows that there is some vertex $v$ of $Z$ such that $g_v$ is
uniquely maximized over $Q$ at $u$, and therefore $u=z_v$.
Suppose that the linear integer programming oracle replied with
an optimal solution to each query. Since $Z$ refines $Q$, this implies
that $Q$ is bounded hence a polytope. Since $c(w_1x,\dots,w_dx)$ is
convex on $\R^n$ and $c$ is convex on $\R^d$, we have that
$$\hskip-3.5cm\max\{c(w_1 x,\dots,w_d x): Ax=b,\ x\in\N^n\}
\ =\ \max\{c(w_1 x,\dots,w_d x): x\in P\}$$
$$\hskip2.5cm\ =\ \max\{c(z): z\in Q\}\ = \ \max\{c(u):u\ \mbox{vertex of }Q\}
\ =\ \max\{c(z_v): v\ \mbox{vertex of } Z\}\ .$$
Using the comparison oracle for $c$, identify that $z_v$ achieving
maximum value $c(z_v)$ over all vertices $v$ of $Z$, and output $x_v$
which is the optimal solution to the convex integer programming problem.
\eoproof

\noindent
Recall that solving the convex integer program means that the
algorithm either returns an optimal solution $x\in\N^n$, or asserts that
the problem is infeasible, or asserts that the polyhedron $\{x\in\R^n\,:\, Ax=b\}$ is unbounded in which case the
problem is generally hopeless (see discussion in the introduction).
It may happen, though, that the projection $Q$ of $P$ is bounded even
though $P$ is not: in this case, there is an optimal solution to the
convex integer programming problem, and our algorithm {\em will} find it.

Lemma \ref{EdgeDirections} bares at once useful consequences for systems
whose defining matrix $A$ is totally unimodular, such as network flow problems
and ordinary ($2$-way) transportation problems. For such totally unimodular
systems, the relevant polyhedron $P$ is integer, that is, we have the equality
$$P\ :=\ \conv\{x\in\N^n:Ax=b\}\ \ =\ \ \{x\in\R^n_+:Ax=b\}\ :=\ L\ .$$
This implies the following two useful properties:
first, for any integer vector $b$, a linear integer programming oracle for
$A,b$ is polynomial time realizable by linear programming over $L$;
and second, a set $E$
covering all edge-directions of $P$ is provided by the set of {\em circuits}
of $A$, that is, minimal-support linear dependencies on the columns of $A$,
whose cardinality is bounded above by $n\choose m$. If $m$ grows slowly, say
$m=O(\log n)$, then this bound is sub-exponential and the algorithm underlying
Lemma \ref{EdgeDirections} might provide a good strategy for addressing
the convex integer maximization problem.

Next, we proceed to prove Theorem \ref{Main}.
We need to recall some definitions.
The {\em Graver basis} of an integer matrix $A$, introduced in
\cite{Graver:75}, is a canonical finite set ${\cal G}(A)$ that can be defined
as follows. Let $\L(A):=\{x\in\Z^n:\ Ax=0\}$ be the lattice of integer
linear dependencies on $A$. Define a partial order $\red$ on $\Z^n$
which extends the coordinate-wise order $\leq$ on $\N^n$ as follows: for two
vectors $u,v\in\Z^n$ put $u\red v$ and say that $u$ is {\em conformal}
to $v$ if $|u_i|\leq |v_i|$ and $u_iv_i\geq 0$ for $i=1,\ldots,n$,
that is, $u$ and $v$ lie in the same orthant of $\R^n$ and each component of
$u$ is bounded by the corresponding component of $v$ in absolute value.
The Graver basis of $A$ is then the set ${\cal G}(A)$ of all
$\red$-minimal vectors in $\L(A)\setminus\{0\}$. For instance, if
$A=(1,2,1)$ then ${\cal G}(A)=\pm\{(2,-1,0),(0,-1,2),(1,0,-1),(1,-1,1)\}$.
For more details on Graver bases and the currently fastest procedure
for computing them see \cite{berndbook,Hemmecke:PSP,4ti2}.

It is known that the universal Gr\"obner bases of $A$, namely the
union of all reduced Gr\"obner bases of the toric ideal of the matrix
$A$, contains all edge directions in the integer hulls within the
polytopes $P_b=\{x: Ax=b, x \geq 0\}$ (see Section 5 in \cite{st}. Since
the Graver bases contains this universal one can deduce the following
property (we include a direct proof here):

\bl{GraverEdgeDirections} For
every integer matrix $A\in\Z^{m\times n}$ and every integer vector
$b\in\N^m$, the Graver basis $\G(A)$ of $A$ covers all edge-directions
of the polyhedron $\conv\{x\in\N^n\,:\,Ax=b\}$ defined by $A$ and $b$.
\el \boproof Consider any edge $e$ of $P:=\conv\{x\in\N^n:Ax=b\}$ and
pick two distinct points $u,v\in e\cap\N^n$. Then $g:=u-v$ is in
$\L(A)\setminus\{0\}$ and hence $g$ is a {\em conformal sum} $g=\sum
g^i$ with $g^i\red g$ and $g^i\in\G(A)$ for all $i$. To see this,
recall that $\G(A)$ is the set of $\red$-minimal elements in
$\L(A)\setminus\{0\}$ and note that $\red$ is a well-ordering; if
$g\in\G(A)$, we are done; otherwise there is an $h\in\G(A)$ with
$h\sqsubset g$ in which case, by induction on $\red$, there is a
conformal sum $g-h=\sum g^i$ giving the conformal sum $g=h+\sum g^i$.

Now, we claim that $u-g^i\in P$ for all $i$. To see this, note first
that $g^i\in\G(A)\subset \L(A)$ implies $Ag^i=0$ and hence $A(u-g^i)=Au=b$;
and second, note that $u-g^i\geq 0$: indeed, if $g^i_j\leq 0$ then
$u_j-g^i_j\geq u_j\geq 0$; and if $g^i_j>0$ then $g^i\red g$ implies
$g^i_j\leq g_j$ and therefore $u_j-g^i_j \geq u_j-g_j=v_j\geq 0$.

Now let $w\in\R^n$ be a linear functional uniquely maximized over $P$
at the edge $e$. Then for all $i$, as just proved, $u-g^i\in P$
and hence $wg^i\geq 0$. But $\sum wg^i= wg=wu-wv=0$,
implying that in fact, for all $i$, we have $wg^i=0$ and therefore
$u-g^i\in e$. This implies that each $g^i$ is a direction of the edge
$e$ (in fact, moreover, all $g^i$ are the same, so $g$ is a
multiple of some Graver basis element).
\eoproof

We also need the following two recent results from \cite{DHOW} on n-fold systems.
The first result builds on stabilization of Graver bases established
by Ho\c{s}ten and Sullivant \cite{HS} and Santos and Sturmfels \cite{SS}.

\bp{GraverComputation}
For any fixed $(r+s)\times t$ integer matrix $A$ there is a polynomial time
algorithm that, given any $n$, computes the Graver basis $\G(A^{(n)})$ of
the n-fold matrix $A^{(n)}=({\bf 1}_n\otimes A_1)\oplus(I_n \otimes A_2)$.
\ep
The second result of \cite{DHOW} combines Proposition \ref{GraverComputation}
and the use of the Graver basis for augmentation.
\bp{Linear}
For any fixed $(r+s)\times t$ integer matrix $A$ there is a polynomial
time algorithm that, given $n$ and vectors $w\in\Z^{nt}$ and
$b\in\Z^{r+ns}$, solves the linear n-fold integer programming problem
$$\max\, \{w x:\ A^{(n)}x=b,\ x\in\N^{nt}\}\ .$$
\ep
Combining Lemma \ref{EdgeDirections}, Lemma \ref{GraverEdgeDirections},
and Propositions \ref{GraverComputation} and \ref{Linear}, we can now
prove Theorem \ref{Main}.

\vskip.2cm\noindent{\bf Theorem \ref{Main}}
{\em For any fixed positive integer $d$ and fixed $(r+s)\times t$ integer
matrix $A$ there is a polynomial\break oracle-time algorithm that,
given $n$, vectors $w_1,\dots,w_d\in\Z^{nt}$ and $b\in\Z^{r+ns}$,
and convex function $c:\R^d\longrightarrow\R$ presented by a
comparison oracle, solves the convex n-fold integer maximization problem
$$\max\, \{c(w_1 x,\dots,w_d x):\ A^{(n)}x=b,\ x\in\N^{nt}\}\ .$$}
\vskip.3cm
\boproof
The algorithm underlying Proposition \ref{Linear} provides a polynomial time
realization of a linear integer programming oracle for $A^{(n)}$ and $b$.
The algorithm underlying Proposition \ref{GraverComputation} allows
to compute the Graver basis $\G(A^{(n)})$ in time which is polynomial
in the input. By Lemma \ref{GraverEdgeDirections},
this set $E:=\G(A^{(n)})$ covers all edge-directions of the polyhedron
$\conv\{x\in\N^{nt}\,:\,A^{(n)}x=b\}$ underlying the convex integer program.
Thus, the hypothesis of Lemma \ref{EdgeDirections} is satisfied and hence
the algorithm underlying Lemma \ref{EdgeDirections} can be used to
solve the convex integer maximization problem in polynomial time.
\eoproof

\section{Applications}
\label{Applications}

We now discuss various applications of our results to multiway
transportation problems, packing problems, vector partitioning and
clustering, extending and unifying applications from \cite{DHOW,HOR,OR,OS}.

\subsection{Multiway transportation problems}
\label{Applications1}

A $k$-way transportation polytope is the set of all $m_1\times\cdots\times m_k$
nonnegative arrays $x=(x_{i_1,\dots,i_k})$ such that the sums of the entries
over some of their lower dimensional sub-arrays (margins) are specified.
More precisely, for any tuple $(i_1,\dots,i_k)$ with
$i_j\in\{1,\dots,m_j\}\cup\{+\}$, the corresponding {\em margin}
$x_{i_1,\dots,i_k}$ is the sum of entries of $x$ over all coordinates
$j$ with $i_j=+$. The {\em support} of $(i_1,\dots,i_k)$ and of
$x_{i_1,\dots,i_k}$ is the set $\supp(i_1,\dots,i_k):=\{j:i_j\neq +\}$
of non-summed coordinates. For instance, if $x$ is a $4\times5\times3\times2$
array then it has $12$ margins with support $F=\{1,3\}$ such as
$x_{3,+,2,+}=\sum_{i_2=1}^5\sum_{i_4=1}^2 x_{3,i_2,2,i_4}$.
Given a family $\F$ of subsets of $\{1,\dots,k\}$ and margin values
$u_{i_1,\dots,i_k}$ for all tuples with support in $\F$,
the corresponding $k$-way transportation polytope is the set of
nonnegative arrays with these margins,
$$T_{\F} \ =\  \left\{\,x\in\R_+^{m_1\times\cdots \times m_k}
\ : \ x_{i_1,\dots,i_k}\,=\,u_{i_1,\dots, i_k}
\,,\ \ \supp(i_1,\dots,i_k)\in\F\,\right\}\ .$$
Transportation polytopes and their integer points (called
contingency tables by statisticians), have been studied and
used extensively in the operations research literature and in
the context of secure statistical data disclosure by public agencies,
see \cite{AT,BR,Cox1,Cox2,DFKPR,QS,Vla,YKK} and references therein.

We now show that when two sides $p$,$q$ of a $3$-way transportation
problem are fixed and one side $n$ is variable, the problem is an
n-fold integer programming problem, and we could therefore conclude
that the convex line-sum $3$-way integer transportation problem is
solvable in polynomial time. Consider the n-fold programming
equations as described after Theorem \ref{Main} in the introduction.
Re-index the arrays as $x=(x^1,\dots,x^n)$ with each
$x^k:=(x^k_{i,j}):=(x_{1,1,k},\dots,x_{p,q,k})$ suitably indexed as a
$pq$ vector representing the $k$-th layer of $x$. Let $r:=t:=pq$ and
$s:=p+q$, and let $A$ be the $(r+s)\times t$ matrix with $A_1:=I_{pq}$
the $pq\times pq$ identity and with $A_2$ the $(p+q)\times pq$ matrix
of equations of the usual $2$-way transportation problem for $p\times
q$ arrays. Finally, define the right-hand side $b=(b^0,b^1,\dots,b^n)$
from the given line-sums by $b^0:=(u_{i,j})$ and
$b^k:=((v_{i,k}),(z_{j,k}))$ for $k=1,\dots,n$. Then the equations
$A_1(\sum_{k=1}^n x^k)=b^0$ represent the constraints
$x_{i,j,+}=u_{i,j}$ of all margins with support $\{1,2\}$, where
summation over layers occurs, whereas the equations $A_2x^k=b^k$ for
$k=1,\dots,n$ represent the constraints $x_{i,+,k}=v_{i,k}$ and
$x_{+,j,k}=z_{j,k}$ of all margins with support $\{1,3\}$ or
$\{2,3\}$, where summations are within a single layer at a time.
Thus, generalizing the recent results of \cite{DHOW} for linear
objective functions, we obtain the following remarkable corollary of
Theorem \ref{Main}.

\vskip.2cm\noindent{\bf Corollary \ref{Threeway}}
{\em For any fixed $d,p,q$ there is a polynomial oracle-time algorithm that,
given $n$, arrays $w_1,\dots,w_d\in\Z^{p\times q\times n}$,
$u\in\Z^{p\times q}$, $v\in\Z^{p\times n}$, $z\in\Z^{q\times n}$,
and convex $c:\R^d\longrightarrow\R$ presented by comparison oracle,
solves the convex integer 3-way transportation problem
$$\max\{\,c(w_1 x, \dots, w_d x)
 \ :\ x\in\N^{p\times q\times n}\,,\ \sum_i x_{i,j,k}=z_{j,k}
\,,\ \sum_j x_{i,j,k}=v_{i,k}\,,\ \sum_k x_{i,j,k}=u_{i,j}\,\}\ .$$}

\noindent
As mentioned before, this is in contrast with the case when the dimensions
of two sides of the tables are variable, in which even the linear integer
$3$-way transportation problem is NP-hard, see \cite{DO1,DO2,DO3}.

The following very general extension of Corollary \ref{Threeway} holds as well.
Consider transportation problems of any fixed dimension $k$ for {\em long}
arrays, namely $m_1\times\cdots\times m_{k-1}\times n$ arrays where
$m_1,\dots, m_{k-1}$ are fixed and only the length (number of layers) $n$
is variable. Further, let $\F$ be any family of subsets of $\{1,\dots,k\}$
(the family of supports of fixed margins). Now re-index the arrays as
$x=(x^1,\dots,x^n)$ with each $x^j=(x_{i_1,\dots,i_{k-1},j})$ a
suitably indexed vector representing the $j$-th layer of $x$.
Then this again is a convex n-fold integer programming
problem with an $(r+s)\times t$ defining matrix $A$, with $t:=\prod m_i$,
with $r,s$, $A_1$ and $A_2$ suitably determined from $\F$, and with the
right-hand side determined from the given margins, in such a way that the
equations $A_1(\sum_{j=1}^n x^j)=b^0$ represent the constraints of all margins
$x_{i_1,\dots,i_k}$ with $i_k=+$ (where summation over layers occurs),
whereas the equations $A_2x^j=b^j$ for $j=1,\dots,n$ represent
the constraints of all margins $x_{i_1,\dots,i_k}$ with $i_k\neq +$
(where summations are within a single layer at a time).
We obtain the following corollary of Theorem \ref{Main} providing
the polynomial time solvability of a very broad class of convex
integer multiway transportation problems.

\bc{Multiway}
For any fixed $d,k,m_1,\dots,m_{k-1}$, and family $\F$ of subsets of
$\{1,\dots,k\}$, there is a polynomial oracle-time algorithm that, given $n$,
arrays $w_1,\dots,w_d\in\Z^{m_1\times \cdots \times m_{k-1}\times n}$,
margin values $u_{i_1,\dots,i_k}$ for all tuples $(i_1,\dots,i_k)$ with
support in $\F$, and convex $c:\R^d\longrightarrow\R$ presented by comparison
oracle, solves the corresponding convex integer multiway transportation problem
$$\max\{\,c(w_1 x, \dots, w_d x)
 \ :\ x\in\N^{m_1\times\cdots \times m_{k-1}\times n}
\,,\ \ x_{i_1,\dots,i_k}\,=\,u_{i_1,\dots, i_k}
\,,\ \ \supp(i_1,\dots,i_k)\in\F\,\}\ .$$
\ec

\subsection{Packing problems}
\label{Applications2}

We consider the following rather general packing problem, which concerns
maximum utility packing of many items of several types in various bins
subject to weight constraints. More precisely, the data is as follows.
There are $t$ types of items. The weight of each item of type $j$ is $v_j$
and there are $n_j$ items of type $j$ to be packed. There are $n$ bins,
where bin $k$ has maximum weight capacity $u_k$. In the linear
version of the problem, there is one utility matrix $w\in\Z^{t\times n}$
where $w_{j,k}$ is the utility of packing one item of type $j$ in bin $k$,
and the objective is to find a feasible packing of maximum total utility.
In the more general convex version, there are $d$ utility matrices
$w_1,\dots,w_d\in\Z^{t\times n}$, representing the packing utilities under
$d$ different criteria. The total utility is the ``balancing" of these
linear utilities under a given convex functional $c$ on $\R^d$.
By incrementing the number $t$ of types by $1$ and suitably augmenting
the data, we may assume that the last type $t$ represents ``slack items"
which occupy the unused capacity in each bin, where the weight of each
slack item is $1$, the utility under each of the $d$ criteria of packing
any slack item in any bin is $0$, and the number of slack bins is the total
residual weight capacity $n_t:=\sum_{k=1}^n u_k-\sum_{j=1}^{t-1} n_jv_j$.
Let $x\in\N^{t\times n}$ be a variable matrix where $x_{j,k}$ represents
the number of items of type $j$ to be packed in bin $k$. Then the convex
packing problem is:
$$\max\{\,c(w_1 x, \dots, w_d x) \ :\ x\in\N^{t\times n}
\,,\ \sum_j v_jx_{j,k}=u_k\,,\ \sum_k x_{j,k}=n_j\,\}\ .$$
By suitably arranging the variables in a vector, it is not hard to
see that this is a convex n-fold integer programming problem with a
$(t+1)\times t$ defining matrix $A$, where $A_1:=I_t$ is the $t\times t$
identity matrix and $A_2:=(v_1,\dots,v_t)$ is a $1\times t$ matrix.
Thus, we obtain the following corollary to Theorem \ref{Main}.

\bc{Packing}
For any fixed number $t$ of types and type weights $v_1,\dots,v_t$,
there is a polynomial oracle-time algorithm that, given $n$, item numbers
$n_j$, bin capacities $u_k$, utilities $w_1,\dots,w_d\in\Z^{t\times n}$,
and convex $c:\R^d\longrightarrow\R$ presented by comparison oracle,
solves the convex integer bin packing problem.
\ec
Note that an interesting special case of bin packing is the classical
{\em cutting stock} problem, and a similar corollary
regarding the solvability of a suitable convex cutting stock problem
can be obtained as well.

\subsection{Vector partitioning and clustering}
\label{Applications3}

The vector partition problem concerns the partitioning of $n$
items among $p$ players to maximize social value subject to
constraints on the number of items each player can receive.
More precisely, the data is as follows. With each item $i$ is
associated a vector $v_i\in\Z^k$ representing its utility under $k$ criteria.
The utility of player $h$ under partition $\pi=(\pi_1,\dots,\pi_p)$ of the
set of items $\{1,\dots,n\}$ is the sum $v^{\pi}_h:=\sum_{i\in\pi_h} v_i$
of utility vectors of items assigned to $h$ under $\pi$.
The social value of $\pi$ is the balancing
$c(v^{\pi}_{1,1},\dots,v^{\pi}_{1,k},\dots,v^{\pi}_{p,1},\dots,v^{\pi}_{p,k})$
of the player utilities, where $c$ is a convex functional on $\R^{pk}$.
In the constrained version, the number $|\pi_h|$ of items that
player $h$ gets is required to be a given number $\lambda_h$
(so $\sum\lambda_h=n$). In the unconstrained version,
there is no restriction on the number of items per player.

Vector partition problems have applications in diverse fields such
as clustering, inventory, reliability, and more - see
\cite{BHR, BH, FOR, HOR, HR, OS, PRW} and references therein.
Here is a typical example.

\be{Clustering}{\bf Minimal variance clustering.}
This is the following problem, which has numerous applications in
the analysis of statistical data: given $n$ observed points
$v_1,\dots,v_n$ in $k$-space, group the points into $p$ clusters
$\pi_1,\dots,\pi_p$ so as to  minimize the sum of cluster variances given by
$$\sum_{h=1}^p{1\over|\pi_h|}\sum_{i\in \pi_h}
||v_i-({1\over|\pi_h|}\sum_{i\in\pi_h} v_i)||^2\ .$$
Consider the instance where there are $n=pm$ points and the
desired clustering is balanced, that is, the clusters should
have equal size $m$. Suitable manipulation of the sum of variances
shows that the problem is equivalent to a constrained partition problem,
where $\lambda_h=m$ for all $h$, and where the convex functional
$c:\R^{pk}\longrightarrow\R$ (to be maximized) is the
Euclidean norm squared, given by
$$c(z)\ =\ ||z||^2\ =\ \sum_{h=1}^p\sum_{i=1}^k|z_{h,i}|^2\ .$$
\ee

If either the number of criteria $k$ or the number of players $p$ is variable,
the partition problem is intractable since it instantly captures NP-hard
problems \cite{HOR}. When both $k,p$ are fixed, both the constrained and
unconstrained versions of the vector partition problem are polynomial time
solvable \cite{HOR,OS}. We now demonstrate how to get this result as a
corollary of Theorem \ref{Main} by showing that both versions are special
convex n-fold integer programming problems. There is an obvious one-to-one
correspondence between partitions and matrices $x\in\{0,1\}^{p\times n}$
with all column-sums equal to one, where partition $\pi$ corresponds to the
matrix $x$ with $x_{h,i}=1$ if $i\in\pi_h$ and $x_{h,i}=0$ otherwise.
Let $d:=pk$ and define $d$ matrices $w_{h,j}\in\Z^{p\times n}$
by setting $(w_{h,j})_{h,i}:=v_{i,j}$ for all $h=1,\dots,p$,
$i=1,\dots,n$ and $j=1,\dots,k$, and setting all other entries to zero.
Then for any partition $\pi$ and its corresponding matrix $x$ we have
$v^{\pi}_{h,j}=w_{h,j}x$ for all $h=1,\dots,p$ and $j=1,\dots,k$.
Therefore, we obtain that the unconstrained vector partition
problem is the convex integer programming problem
$$\max\{\,c(w_{1,1} x, \dots, w_{p,k} x) \ :\ x\in\N^{p\times n}
\,,\ \sum_h x_{h,i}=1\,\}\ .$$
Suitably arranging the variables in a vector, it is not hard to
see that this is a convex n-fold integer programming problem
with a $(0+1)\times p$ defining matrix $A$, where $A_1$ is empty
and $A_2:=(1,\dots,1)$. Similarly, the constrained vector partition
problem is the convex integer programming problem
$$\max\{\,c(w_{1,1} x, \dots, w_{p,k} x) \ :\ x\in\N^{p\times n}
\,,\ \sum_h x_{h,i}=1\,,\ \sum_i x_{h,i}=\lambda_h\,\}\ .$$
Again, it can be seen that this is a convex n-fold integer programming
problem, now with a $(p+1)\times p$ defining matrix $A$, where now $A_1:=I_p$
is the $p\times p$ identity matrix, and $A_2:=(1,\dots,1)$ as before.

Thus, we obtain the following corollary to Theorem \ref{Main}.
\bc{Partition}
For any fixed number $p$ of players and number $k$ of criteria,
there is a polynomial oracle-time algorithm that, given $n$, item vectors
$v_i\in\Z^k$, positive integers $\lambda_h$,
and convex $c:\R^{pk}\longrightarrow\R$ presented by comparison oracle,
solves the constrained and the unconstrained vector partition problems.
\ec

\section*{Acknowledgement}

We thank an anonymous referee for helpful suggestions
that improved the presentation of the paper.

\noindent {\small Jesus De Loera}\newline
\emph{University of California at Davis, Davis, CA 95616, USA}\newline
\emph{email: deloera{\small @}math.ucdavis.edu},
\ \ \emph{http://www.math.ucdavis.edu/{\small $\sim$deloera}}

\noindent {\small Raymond Hemmecke}\newline
\emph{Otto-von-Guericke Universit\"at Magdeburg,
D-39106 Magdeburg, Germany}\newline
\emph{email: hemmecke{\small @}imo.math.uni-magdeburg.de},
\ \ \emph{http://www.math.uni-magdeburg.de/{\small $\sim$hemmecke}}

\noindent {\small Shmuel Onn}\newline
\emph{Technion - Israel Institute of Technology, 32000 Haifa, Israel}\newline
\emph{email: onn{\small @}ie.technion.ac.il},
\ \ \emph{http://ie.technion.ac.il/{\small $\sim$onn}}

\noindent {\small Uriel G. Rothblum}\newline
\emph{Technion - Israel Institute of Technology, 32000 Haifa, Israel}\newline
\emph{email: rothblum{\small @}ie.technion.ac.il},
\ \ \emph{http://ie.technion.ac.il/{\small rothblum.phtml}}

\noindent {\small Robert Weismantel}\newline
\emph{Otto-von-Guericke Universit\"at Magdeburg,
D-39106 Magdeburg, Germany}\newline
\emph{email: weismantel{\small @}imo.math.uni-magdeburg.de},
\ \ \emph{http://www.math.uni-magdeburg.de/{\small $\sim$weismant}}

\end{document}